\documentclass[11pt]{article}

\usepackage{amssymb,amsmath, amsfonts, amsthm}

\usepackage{graphicx}
\usepackage{pgfplots}
\usepackage{tikz}
\usepackage{caption}
\usepackage{booktabs}
\usepackage{array}
\usepackage{verbatim}
\usepackage{subfig}
\usepackage{geometry}
\usepackage{fancyhdr}

\newcommand{\vertex}{\node[vertex]}
\tikzstyle{vertex}=[circle, draw, inner sep=0pt, minimum size=6pt]
\newtheorem{theorem}{Theorem}
\newtheorem{lemma}{Lemma}
\newtheorem{corollary}{Corollary}
\newtheorem{obs}{Observation}
\newtheorem{prop}{Proposition}

\newcommand{\smallqed}{{\tiny ($\Box$)}}

\begin{document}

\title{Well-edge-dominated graphs containing triangles}
\author{
Jake Berg$^{b}$\and   Perryn Chang$^{b}$  \and Claire Kaneshiro$^{c}$ \and
Kirsti Kuenzel$^{a}$ \and Ryan Pellico$^{a}$ \and Isabel Renteria$^{d}$\and Sumi Vora$^{e}$ }

\maketitle

\begin{center}
$^a$ Department of Mathematics, Trinity College, Hartford, CT\\
$^b$ Department of Mathematics, Yale University, New Haven, CT\\
$^c$ Department of Mathematics, Princeton University, Princeton, NJ\\
$^d$ Department of Mathematics and Statistics,  Loyola University Chicago, Chicago, IL\\
$^e$ Department of Mathematics and Statistics,  Pomona College, Claremont, CA\\
\end{center}
\medskip

\maketitle
\begin{abstract}
A set of edges $F$ in a graph $G$ is an edge dominating set if every edge in $G$ is either in $F$ or shares a vertex with an edge in $F$. $G$ is said to be well-edge-dominated if all of its minimal edge dominating sets have the same cardinality. Recently it was shown that any triangle-free well-edge-dominated graph is either bipartite or in the set $\{C_5, C_7, C_7^*\}$ where $C_7^*$ is obtained from $C_7$ by adding a chord between any pair of vertices distance three apart. In this paper, we completely characterize all well-edge-dominated graphs containing exactly one triangle, of which there are two infinite families. We also prove that there are only eight well-edge-dominated outerplanar graphs, most of which contain at most one triangle. 
\end{abstract}

\maketitle
\section{Introduction}
Given a graph $G$,  a set of edges $F$ in $G$ is called a \emph{matching} if no pair of edges in $F$ share a common vertex. If every edge of $G$ is either in $F$ or adjacent to an edge in $F$, then $F$ is an \emph{edge dominating set} of $G$. A graph is said to be \emph{equimatchable} if all of of its maximal matchings have the same size, and is said to be \emph{well-edge-dominated} if all of its minimal edge dominating sets have the same size. Equimatchable graphs and well-edge-dominated graphs are the edge version counterparts to well-covered graphs and well-dominated graphs. In this paper, we are interested in studying well-edge-dominated graphs, which as a set is properly contained in the set of all equimatchable graphs. Equimatchable graphs were first studied independently by Lewin \cite{L-1974} and Meng \cite{M-1974}, and later Lesk et al. \cite{LPP-1984}. Frendrup, Hartnell, and Vestergaard \cite{FHV-2010} proved that a connected equimatchable graph of girth at least $5$ is either $C_5$, $C_7$ or in a specific class of bipartite graphs. They also proved that any equimatchable graph with girth at least $5$ is also well-edge-dominated. Therefore, all well-edge-dominated graphs of girth $5$ or more are completely characterized. In 2023,  B\"{u}y\"{u}k\c{c}olak et al. \cite{BGO-2023} provided a characterization of the connected, triangle-free equimatchable graphs that are not bipartite. Anderson et al. \cite{AKR-2022} identified which of those triangle-free equimatchable graphs are in fact well-edge-dominated. Therefore, to finish the characterization of all non-bipartite well-edge-dominated graphs, one only needs to focus on those containing triangles. 

As a first step in this direction, we define two infinite families of graphs, both of which are obtained from a bipartite well-edge-dominated graph. We use the following terminology to define these families. Let $G$ be a well-edge-dominated bipartite graph with bipartition $V(G) = A \cup B$ where $|A| < |B|$. We say that $w \in B$ is \emph{detachable} if $G - w$ is also well-edge-dominated. Additionally, we say $w$ is \emph{strongly detachable} if it is detachable and each vertex in $N_G(w)$ is a support vertex (adjacent to a vertex with degree $1$) in $G-w$. A graph $G$ is in the family $\mathcal{T}$ if it is obtained from the disjoint union of $K_3$ and a well-edge-dominated bipartite graph $G' = (A\cup B, E)$ where $|A| < |B|$ by identifying $z \in V(K_3)$ with $w \in V(G')$ where $w$ is detachable. A graph $G$ is in the family $\mathcal{F}$ if it is obtained from the disjoint union of the house graph $\mathcal{H}$, depicted in Figure~\ref{fig:houses}, and a well-edge-dominated bipartite graph $G' = (A \cup B, E)$ where $|A| < |B|$ by identifying the vertex of degree $2$ on the triangle in $\mathcal{H}$ with $w \in V(G')$ where $w$ is strongly detachable. The majority of this paper is spent proving the following (where $\{Cr, \mathcal{H}, \mathcal{DH}\}$ are depicted in Figure~\ref{fig:houses}). 
\begin{theorem} \label{thm:onetriangleWED}
 $G$ is a connected, well-edge-dominated graph with exactly one triangle if and only if $G \in \mathcal{T}\cup \mathcal{F} \cup \{K_3, Cr, \mathcal{H}, \mathcal{DH}\}$.
\end{theorem}

The remainder of the paper is organized as follows. In Section~\ref{sec:defn}, we provide useful definitions and previous results pertinent to this paper. In Section~\ref{sec:one}, we provide the proof of Theorem~\ref{thm:onetriangleWED}. In Section~\ref{sec:OP} we characterize the well-edge-dominated outerplanar graphs. 

  \subsection{Definitions and Previous Results}\label{sec:defn}
  Throughout this paper, we consider only simple, finite, undirected graphs. Let $G = (V(G), E(G))$ be any graph. We let $n(G)$ denote the order of the graph, namely $|V(G)|$. We say that $G$ is \emph{nontrivial} if $n(G) \ge 2$. We use the notation $[k] = \{1, 2, \dots, k\}$ for any positive integer $k$. Given a vertex $v \in V(G)$, the \emph{neighborhood of $v$} in $G$ is the set of vertices adjacent to $v$, and denoted $N_G(v)$. The \emph{closed neighborhood} of $v$ is defined as $N_G[v] = N_G(v) \cup \{v\}$ and the degree of $v$ is precisely $\deg_G(v) = |N_G(v)|$. A \emph{leaf} is a vertex of degree $1$ and a \emph{support vertex} is any vertex adjacent to a leaf. Further, $v$ is a \emph{strong support vertex} if $v$ is adjacent to at least two leaves. The length of the smallest cycle in $G$ is referred to as the \emph{girth} of $G$ and denoted $g(G)$. We say that $x \in V(G)$ is a \emph{cut-vertex} if $G$ is connected and $G-x$ is disconnected. 
  
  Given any pair of edges $e$ and $f$ in $E(G)$, we say $e$ and $f$ are \emph{adjacent} if they share a common vertex. The \emph{closed edge neighborhood} of $f$ is the set $N_e[f]$ consisting of $f$ together with all edges in $G$ that are adjacent to $f$. Similarly, for any $F \subseteq E(G)$, the \emph{closed edge neighborhood} of $F$ is the set $N_e[F]$ defined by $N_e[F]=\cup_{f\in F}N_e[f]$. The edge $f$ is said to \emph{dominate} the set $N_e[f]$. An edge $g$ is called a \emph{private edge neighbor} of $f$ with respect to $F$ if $g \in N_e[f]$ and $g \not\in N_e[F-\{f\}]$. If $N_e[F]= E(G)$, then $F$ is called an \emph{edge dominating set} of $G$.  The \emph{edge domination number} of $G$, denoted by $\gamma_e(G)$, is the smallest cardinality of an edge dominating set in $G$. A \emph{matching} in $G$ is a set of independent edges, i.e., a set in which no pair of edges share a common vertex. The maximum cardinality among all matchings in $G$ is the matching number of $G$, denoted $\alpha'(G)$. Given a set of edges $M\subseteq E(G)$, any vertex incident to an edge of $M$ is said to be \emph{saturated} by $M$. Similarly, if $U$ is any set of vertices, we say that $M$ saturates $U$ if each vertex of $U$ is saturated by $M$. A \emph{perfect matching} is a matching which saturates all vertices in $G$. 
  
 A graph $G$ is called \emph{equimatchable} if all of its maximal matchings in $G$ have the same cardinality. Additionally, $G$ is well-edge-dominated if all of its minimal edge dominating sets have the same cardinality. As any matching in $G$ is a minimal edge dominating set, if $G$ is well-edge-dominated, then it is also equimatchable.

 Frendrup et al. pointed out in \cite{FHV-2010} that every  equimatchable graph of girth $5$ or more is also well-edge-dominated, and they proved the following.
 
 \begin{theorem}\cite{FHV-2010}\label{thm:girth5} If $G$ is a connected graph with $g(G) \ge 5$, then $G$ is well-edge-dominated if and only if $G \in \{K_2, C_5, C_7\}$ or $G$ is bipartite with partite sets $U$ and $V$ such that $U$ is the set of all support vertices of $G$. 
 \end{theorem}
 
 Thus, the above characterizes all well-edge-dominated graphs with girth $5$ or more. Then Anderson et al. \cite{AKR-2022} showed the following (where $C_7^*$ is the graph obtained from $C_7$ by adding a chord between two vertices distance three apart). 
 
  \begin{theorem}\label{thm:girth5red}\cite{AKR-2022} If $G$ is a connected, well-edge-dominated graph with $g(G) \ge 4$, then either $G$ is bipartite or $G \in \{C_5, C_7, C_7^*\}$.
  \end{theorem}
  
Equimatchable bipartite graphs containing $4$-cycles have been studied in \cite{BGO-2020}. Notably, the following is shown. 

  \begin{lemma}\cite{BGO-2020} Let $G = (U\cup V, E)$ with $|U| < |V|$ be a connected equimatchable bipartite graph. Then each vertex $u \in U$ satisfies at least one of the following.
  \begin{enumerate}
  \item[(i)] $u$ is a support vertex in $G$, or
  \item[(ii)] $u$ is on a $4$-cycle.
  \end{enumerate}
  \end{lemma}

  Note that the above assumes that $G= (U \cup V, E)$ is an equimatchable  bipartite graph with $|U| < |V|$. Sumner \cite{S-1979} had already characterized all equimatchable bipartite graphs with $|U| = |V|$.  A graph $G$ is defined to be \emph{randomly matchable} if it is an equimatchable graph admitting a perfect matching. Sumner, in fact,  characterized all randomly matchable graphs. 
  
  \begin{theorem}\cite{S-1979} A connected graph is randomly matchable if and only if it is isomorphic to $K_{2n}$ or $K_{n,n}$ for $n \ge 1$. 
  \end{theorem}
  We also have the following result from \cite{DE-2019}. 
    \begin{theorem}\cite{DE-2019} Let $G = (U \cup V,E)$ be a connected bipartite graph with $|U| \le |V|$. Then $G$ is equimatchable if and only if every maximal matching of $G$ saturates all vertices in $U$.
  \end{theorem}
  
  Therefore, if $G$ is an equimatchable bipartite graph $G = (U\cup V, E)$ where $|U| = |V|$, then it contains a maximal matching that saturates $U$ and therefore $V$, so it admits a perfect matching. This implies that $G = K_{n,n}$ for $n \ge 1$. Since every well-edge-dominated graph is equimatchable, we have the following corollary.
  
  \begin{corollary}\label{cor:Knn} If $G$ is a well-edge-dominated bipartite graph $G = (U\cup V, E)$ where $|U| = |V|$, then $G = K_{n,n}$ for $n \ge 1$. 
  \end{corollary}
  
  Throughout the paper, we will repeatedly use the following facts. 
  
  \begin{lemma}\cite{AKR-2022}\label{lem:reduce} Let $M$ be any matching in a graph $G$. If $G$ is well-edge-dominated, then $G - N_e[M]$ is well-edge-dominated. If $G$ is equimatchable, then $G - N_e[M]$ is equimatchable. 
  \end{lemma}

  We also observe the following about well-edge-dominated bipartite graphs.

  \begin{obs}\label{obs:matchings} Let $G= (U\cup V, E)$ be an equimatchable bipartite graph where $|U| <|V|$. For any $v \in V$, there exist two maximal matchings $M_1$ and $M_2$ in $G$ such that $M_1$ saturates $v$ and $M_2$ does not saturate $v$. 
  \end{obs}
  
  We also use the following result first shown in \cite{FHV-2010}.
  
\begin{lemma}\cite{FHV-2010}\label{lem:exmatch} Let $G$ be an equimatchable graph. If $M_1$ and $M_2$ are matchings in $G$ and $N_e[M_1] = N_e[M_2]$, then $|M_1| = |M_2|$.
\end{lemma}
\section{Graphs containing exactly one triangle}\label{sec:one}
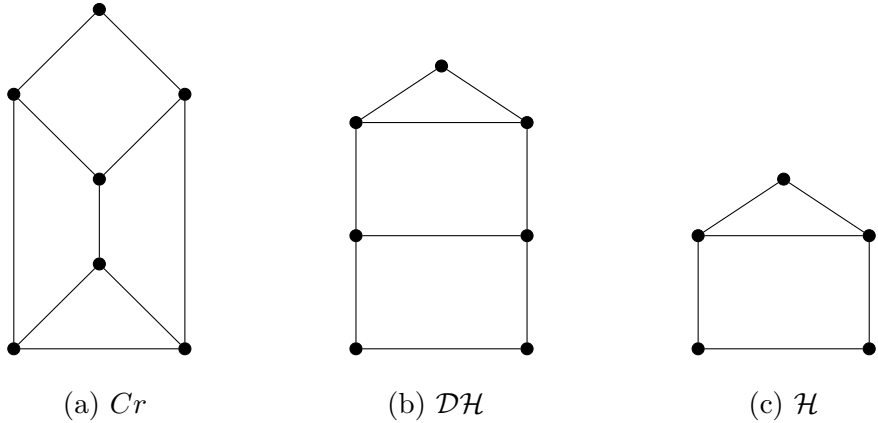
\begin{figure}[h!]
\begin{center}
\begin{tikzpicture}[scale=.75]
    \vertex (1) at (0,0)  [fill, scale=.75]{};
    \vertex (2) at (0,2) [fill, scale=.75]{};
    \vertex (3) at (3,0)  [fill, scale=.75]{};
    \vertex (4) at (3,2)   [fill, scale=.75]{};
    \vertex (5) at (0,4) [fill, scale=.75]{};
    \vertex (6) at (3,4)  [fill, scale=.75]{};
    \vertex (7) at (1.5, 5)  [fill, scale=.75]{};
    \vertex (8) at (6,0) [fill, scale=.75]{};
    \vertex (9) at (9,0)  [fill, scale=.75]{};
    \vertex (10) at (6,2)  [fill, scale=.75]{};
    \vertex (11) at (9,2)  [fill, scale=.75]{};
    \vertex (12) at (7.5,3)  [fill, scale=.75]{};
      \vertex (1b) at (-6, 0)  [scale=.75, fill=black]{};
    \vertex (2b) at (-4.5,1.5)  [scale=.75, fill=black]{};
    \vertex (3b) at (-3,0)  [scale=.75, fill=black]{};
    \vertex (4b) at (-4.5, 3)  [scale=.75, fill=black]{};
    \vertex (5b) at (-6, 4.5)  [scale=.75, fill=black]{};
    \vertex (6b) at (-4.5, 6)  [scale=.75, fill=black]{};
    \vertex (7b) at (-3, 4.5)  [scale=.75, fill=black]{};
    \node(A) at (1.4, -1)[]{(b) $\mathcal{DH}$};
    \node(B) at (7.5, -1)[]{(c) $\mathcal{H}$};
    \node(C) at (-4.4, -1)[]{(a) $Cr$};

    \path 
    (1) edge (2)
    (1) edge (3)
    (4) edge (2)
    (4) edge (3)
    (5) edge (2)
    (6) edge (4)
    (7) edge (6)
    (7) edge (5)
    (5) edge (6)

    (8) edge (9)
    (8) edge (10)
    (11) edge (10)
    (11) edge (9)
    (11) edge (12)
    (12) edge (10)
    
        (1b) edge (2b)
    (2b) edge (3b)
    (3b) edge (1b)
    (4b) edge (2b)
    (4b) edge (5b)
    (5b) edge (6b)
    (4b) edge (7b)
    (7b) edge (3b)
    (1b) edge (5b)
    (6b) edge (7b)
    ;

\end{tikzpicture}
\caption{The crystal graph, the dream house, and the house graphs}
\label{fig:houses}
\end{center}
\end{figure}

To prove Theorem 1, we first define two infinite families of graphs as well as a few exceptional graphs. Define $\mathcal{H}$, $\mathcal{DH}$, and $Cr$ to be the three graphs depicted in Figure~\ref{fig:houses}. $\mathcal{H}$ is the usual house graph, we refer to $\mathcal{DH}$ as the ``dream house", and we refer to $Cr$ as the ``crystal" graph.

To define the infinite families, we rely heavily on well-edge-dominated bipartite graphs. We will use the following terminology throughout this section. Let $G$ be  a well-edge-dominated bipartite graph with bipartition $V(G) = A \cup B$ where  $|A| < |B|$. We say that $w \in B$ is \emph{detachable} if $G-w$ is well-edge-dominated. Note that in this case $\gamma_e(G-w) = |A|$ as $G-w$ contains a matching of size $|A|$. Additionally, we say that $w \in B$ is \emph{strongly detachable} if $G-w$ is well-edge-dominated and each vertex in $N_G(w)$ is a support vertex in $G-w$.

 The first family $\mathcal{F}$ is defined as follows. We say that $G \in \mathcal{F}$ if $G$ is obtained from the disjoint union of the house graph $\mathcal{H}$ and a well-edge-dominated bipartite graph $G'$ with bipartition $ A \cup B$,   $|A|< |B|$,  by identifying the vertex of degree $2$  on the triangle in $\mathcal{H}$ with a strongly detachable vertex $w \in V(G')$. 

The second family $\mathcal{T}$ is defined as follows. We say that $G \in \mathcal{T}$ if $G$ is obtained from the disjoint union of $K_3$ and a well-edge-dominated bipartite graph $G'$ with bipartition  $A \cup B$,  $|A| < |B|$,  by identifying $z \in V(K_3)$ with $w \in V(G')$ where $w$ is detachable. We first show that every graph in $\mathcal{T} \cup \mathcal{F}$ is well-edge-dominated.

\begin{figure}
\centering
\begin{tikzpicture}[scale=2]
    \vertex (1) at (0,0)  [scale=.75, fill]{};
    \vertex (2) at (0.5,0.5) [scale=.75, fill]{};
    \vertex (3) at (1,0)  [scale=.75, fill]{};
    \vertex (4) at (1.5,0.5)  [scale=.75, fill]{};
    \vertex (5) at (1.5,0)   [scale=.75, fill]{};
    \vertex (6) at (2,0)   [scale=.75, fill]{};
    \vertex (7) at (1,-.7)   [scale=.75, fill]{};
    \vertex (8) at (2,-.7)  [scale=.75, fill]{};

    \path 
    (1) edge (2)
    (2) edge (3)
    (3) edge (4)
    (4) edge (5)
    (4) edge (6)
    (7) edge (8)
    (8) edge (5)
    (7) edge (5)
    ;

\end{tikzpicture}
\caption{Example of a graph in $\mathcal{T}$}
\end{figure}
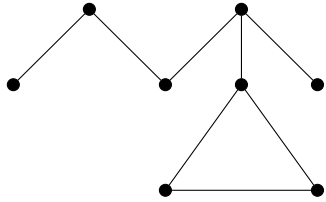

Now we are ready to prove one direction of Theorem 1. 

\begin{prop} If $G \in \mathcal{T} \cup \mathcal{F}\cup \{K_3, Cr, \mathcal{H}, \mathcal{DH}\}$, then $G$ is well-edge-dominated. 
\end{prop}

\begin{proof}
One can easily verify that $K_3, Cr, \mathcal{H}$, and $\mathcal{DH}$ are well-edge-dominated. Suppose first that $G \in \mathcal{T}$. We shall  assume $G$ is constructed from $K_3 = xyz$ and $G'$ with bipartition $A\cup B$, where $|A|< |B|$, by identifying $z$ and $w \in B$ where $w$ is detachable. We will refer to the contracted vertex in $G$ upon identifying $z$ and $w$ as $z$. We know there exists a maximal matching in $G$ of cardinality $|A|+1$. Let $F$ be any minimal edge dominating set of $G$. It follows that $F$ contains some edge from $xyz$ for otherwise $xy$ is not dominated. Suppose first that $xy \in F$. Note that $xz$ and $yz$ cannot be in $F$ for otherwise $F - \{xy\}$ is also an edge dominating set. Thus, $F - \{xy\}$ is a minimal edge dominating set of $G'$, implying $|F- \{xy\}| = |A|$. So assume that $xy \not\in F$. It follows that either $xz$ or $yz$ is in $F$, but not both. With no loss of generality, we may assume that $xz \in F$. Thus, $F - \{xz\}$ is a minimal edge dominating set of $G'-w$, which is assumed to be well-edge-dominated with $\gamma_e(G'-w) = |A|$. Thus, $|F - \{xz\}| = |A|$. In either case, $|F| = |A|+1$ and $G$ is well-edge-dominated. 

Next, suppose that $G \in \mathcal{F}$. We shall assume $G$ is constructed from $\mathcal{H}$ with the degree $2$ vertex on the triangle being $z$ and bipartite graph $G'$ with bipartition $A\cup B$ such that $|A| < |B|$ and $w \in B$ is strongly detachable. We will refer to the contracted vertex in $G$ upon identifying $z$ and $w$ as $z$. Note first that there exists a maximal matching in $G$ of cardinality $|A|+2$. Next, let $F$ be any minimal edge dominating set of $G$. We know $F$ contains at least two edges from $\mathcal{H}$. Label the remaining vertices in $\mathcal{H}$ as $x, y, t, s$ where $x$ and $y$ are on the triangle and $t$ is adjacent to $x$. 

Suppose first that $F$ contains exactly two edges from $\mathcal{H}$. It must be that at most one of these edges is incident to $z$ for otherwise $st$ is not dominated by $F$. If neither edge is incident to $z$, then we know that $F \cap E(G')$ is a minimal edge dominating set of $G'$, implying $|F \cap E(G')| = |A|$ so $|F| = |A| + 2$. Therefore, we shall assume exactly one edge in $F\cap E(\mathcal{H})$ is incident to $z$. Without loss of generality, we may assume $xz \in F$. It follows that $F \cap E(G')$ is a minimal edge dominating set of both $G'$ and $G' - w$, which satisfies $\gamma_e(G' - w) = |A|$. Thus, $|F| = |A|+2$. 

Finally, assume that $F$ contains exactly three edges from $\mathcal{H}$. The only way this is possible is if all three edges are incident to a common vertex, either $x$ or $y$. Without loss of generality, assume all three edges are incident to $x$. Since $F$ is minimal, $xz$ must have a private neighbor with respect to $F$ of the form $za$ where $a \in A$. However, $a$ is a support vertex in $G' - w$ by assumption meaning that $a$ must be saturated by $F$. This contradiction shows that no such minimal edge dominating set $F$ containing three edges from $\mathcal{H}$ exists and $G$ is in fact well-edge-dominated. 

\end{proof}

For the remainder of this section, we prove the other direction of Theorem 1; namely, that the only well-edge-dominated graphs containing exactly one triangle are in $\mathcal{T} \cup \mathcal{F} \cup \{K_3, Cr, \mathcal{H}, \mathcal{DH}\}$. 

\subsection{Proof of Theorem 1}
To begin, we start with a few necessary conditions for a graph containing exactly one triangle to be well-edge-dominated. We begin with a preliminary result.
\begin{lemma}\label{lem:noleaves} If $G$ is a graph of order at least $4$ obtained from the triangle $xyz$ by appending leaves to $x, y$, or $z$, then $G$ is not well-edge-dominated. 
\end{lemma}

\begin{proof}
Suppose to the contrary that $G$ is obtained from $xyz$ by appending at least one leaf, call it $\ell_z$, to $z$. If $\deg_G(x) = 2$, then $\{yz\}$ and $\{xy, z\ell_z\}$ are maximal matchings in $G$, which is a contradiction to the fact that $G$ is equimatchable. Thus, we may assume that $\deg_G(x)\ge 3$, and similarly $\deg_G(y) \ge 3$. Let $\ell_x$ be a leaf of $x$ and let $\ell_y$ be a leaf of $y$. It follows that $\{xy, z\ell_z\}$ and $\{x\ell_x, y\ell_y, z\ell_z\}$ are maximal matchings in $G$, another contradiction.
\end{proof}

The above result implies that for the remainder of this section, if $G$ is a well-edge-dominated graph containing exactly one triangle $xyz$, then the graph induced by $G - \{x, y, z\}$ contains a nontrivial component. We also note that throughout this section, we shall assume that $G$ is connected. Next, we point out the following property of well-edge-dominated bipartite graphs containing a cut-vertex. 

\begin{lemma}\label{lem:subbip2} Let $G = (A \cup B, E)$ be a connected, well-edge-dominated bipartite graph where $|A|<|B|$. If $x \in B$ is a cut-vertex of $G$, then  $G-x$ is well-edge-dominated and we can choose a bipartition $A' \cup B'$ of $G-x$ where $|A'| \le |B'|$ and $A' = A$. 
\end{lemma}

\begin{proof} Note first that $G-x = H_1\cup \cdots \cup H_k$ where each $H_i$ is a component of $G-x$. Pick any edge $xy$ where $y \in V(H_i)$. We know by Lemma~\ref{lem:reduce} that $G-N_e[xy]$ is well-edge-dominated, implying that each $H_j$ where $j \ne i$ is well-edge-dominated. Again choosing some $y \in V(H_j)$ where $j \ne i$, the same argument shows $H_i$ is also well-edge-dominated. It follows that $G-x$ is well-edge-dominated. As each $H_i$ is bipartite, we can partition $V(H_i)$ as $A_i \cup B_i$ where $|A_i| \le |B_i|$. Moreover, if $|A_i| = |B_i|$, then we choose $A_i$ to be the subset of $A$. Suppose for some $i \in [k]$, $B_i \subset A$ and $|A_i| < |B_i|$. Reindexing if necessary, we may assume $i=1$. Let $t$ be a neighbor of $x$ in $H_1$. By Observation~\ref{obs:matchings}, there exist two maximal matchings of $H_1$ with the same cardinality, $M_1$ which saturates $t$ and $M_2$ which does not saturate $t$. Since each component of $G-x$ is well-edge-dominated, we can pick maximal matchings in each $H_i$ for $2 \le i \le k$ which saturate the neighbors of $x$ in $H_i$. Call the resulting matching $F$. It follows that $F \cup M_1$ and $F \cup M_2 \cup \{xt\}$ are maximal matchings in $G$, which is a contradiction. Therefore, the result holds.
\end{proof}

Finally, we prove one last lemma about well-edge-dominated bipartite graphs that is needed before considering well-edge-dominated graphs containing a triangle.

\begin{lemma}\label{lem:subbip} Let $G=(A\cup B, E)$ be a connected, well-edge-dominated bipartite graph where $|A| \le |B|$. If $x$ is a cut-vertex in $G$ where $G - x$ is the disjoint union of components $H_1, \dots, H_k$, then the following hold.
\begin{enumerate}
\item[1.]  If $x \in A$, then for any $i\in[k]$, there exists a maximal matching in $G_i$, the graph induced by $ V(H_i)\cup\{x\}$, which saturates the partite set containing $x$ in $G_i$. 
\item[2.] If $x \in B$, then for any subset $H_{\alpha_1}, \dots, H_{\alpha_j}$ of components of $G-x$ where $j<k$, there exist two maximal matchings $M_1$ and $M_2$ of the same size in the graph induced by $V(G) - \cup_{i=1}^j V(H_{\alpha_i})$ where $M_1$ saturates $x$ and $M_2$ does not saturate $x$. 
\end{enumerate}
\end{lemma}

\begin{proof} Note first that if $|A| = |B|$, then $G = K_{n,n}$ which does not contain a cut-vertex. Therefore, we may assume $|A| < |B|$. 
Suppose first that $x \in A$ and consider the graph induced by $V(H_i) \cup \{x\}$ for some $i\in [k]$, call it $G_i$. Note that $G_i$ is necessarily bipartite. Let $e$ be any edge of $G$ incident to $x$ and a vertex in $V(H_i)$. As $G$ is well-edge-dominated, we can extend $e$ to a maximal matching $M$ in $G$ that saturates $A$. Let $M_{G_i}$ be those edges in $M$ that are also in $G_i$. It follows that $M_{G_i}$ is a maximal matching in $G_i$ which saturates the partite set containing $x$ as there are no edges between a vertex  in $V(H_i)$ and a vertex in $V(G) - V(G_i)$. 

Next, suppose $x \in B$ and consider the graph induced by $V(G) - \left(\cup_{i=1}^j V(H_i)\right)$ where $j<k$, as above call it $J$. Note first that for each $i \in [k]$, we can partition $V(H_i) = C_i \cup D_i$ where each of $C_i$ and $D_i$ are independent sets in $G$. Relabeling if necessary, we shall assume for each $i \in [k]$ that $x$ has a neighbor in $C_i$. Further, $x$ has no neighbor in $D_i$ as this would imply that $G$ contains an odd cycle which cannot be. Now by Lemma~\ref{lem:subbip2},  each $H_i$ is well-edge-dominated and we may assume $|C_i| \le |D_i|$ and $C_i \subset A$ for each $i \in [k]$. As $J$ is bipartite, we can write $V(J) = A' \cup B'$ where $A'\subset A$ and $x \in B'$. Let $e$ be any edge of $G$ incident to $x$ and a vertex in $V(H_{j+1}) \cup \cdots \cup V(H_k)$. Since $G$ is well-edge-dominated, we can extend $e$ to a maximal matching $M_1$ in $G$ which saturates $A$. Let $M_1^J$ be those edges in $M_1$ that are also in $J$. It follows that $M^J_1$ is a maximal matching in $J$ that saturates $A'$. On the other hand, we can choose an edge $f$ of $G$ incident to $x$ and a vertex in $H_1$ and extend $f$ to a maximal matching $M_2$ in $G$ which saturates $A$. If $M^J_2$ are those edges of $M_2$ that are also in $J$, it follows that $M^J_2$ is a maximal matching in $J$ that saturates each vertex in $A'$ but does not saturate $x$. Thus, $|M^J_1| = |A'| = |M^J_2|$. 
\end{proof}

We use Lemma~\ref{lem:subbip} to prove the following result, which is ultimately the backbone of the proof of Theorem 1. 

\begin{theorem}\label{thm:TRI1} If  $G\ne K_3$ is a connected, well-edge-dominated  graph containing exactly one triangle $xyz$, then the component of $G'=G-N_e[xy]$ containing $z$, call it $G'_z$, is bipartite. Moreover, if $G'_z$ is nontrivial, then one of the following is true:
\begin{enumerate}
\item[(a)] every edge incident to $x$ other than $xy$ and $xz$, and every edge incident to $y$ other than $xy$ and $yz$, if such an edge exists,  is incident to a vertex in $G_z'$. Moreover, if $\deg_G(x)= 2$ or $\deg_G(y) = 2$, then $z$ is detachable in $G_z'$.
\item[(b)] every component of $G'$ other than $G_z'$ has a vertex adjacent to $x$ in $G$, $\deg_G(y) \ge3$ and $z$ is detachable in $G_z'$.
\end{enumerate}
\end{theorem}

\begin{proof}
In what follows, we let $G^{xy}=G-N_e[xy]$, $G^{xz} = G - N_e[xz]$, and $G^{yz} = G - N_e[yz]$. For each pair $ab \in \{xy, xz, yz\}$, enumerate the components of $G^{ab}$ as $G^{ab}_1, G_2^{ab}, \dots, G_{k_{ab}}^{ab}$ where we shall assume that $c = \{x, y, z\} - \{a, b\}$ is contained in $G_1^{ab}$. Note that for each $ab \in \{xy, xz, yz\}$, $G^{ab}$  is a well-edged-dominated graph. We first consider $G^{xy}$. 

Since $g(G^{xy}) \ge 4$, we know by Theorem~\ref{thm:girth5red} that each component of $G^{xy}$ is either bipartite or in $\{C_5, C_7, C_7^*\}$. Suppose first that $G^{xy}_1=C_5$.
Label the vertices of $G^{xy}_1$ as $\{v_1, v_2, v_3, v_4,v_5\}$ with $v_1=z$ and $v_iv_{i+1} \in E(G)$ for $i \in[5]$. Consider the two matchings  $N_1=\{xz, v_2v_3, v_4v_5\}$ and $N_2=\{xz, v_3v_4\}$. Note that $\deg_G(v_5) = 2 = \deg_G(v_2)$ as we have assumed that  $G$ contains exactly one triangle (so there are no edges between $\{v_2, v_5\}$ and $\{x, y\}$). Moreover, $N_G(v_4) \subset \{v_3, v_5, x, y\}$ and $N_G(v_3) \subset \{v_2, v_4, x, y\}$. Thus,  $N_e[N_1]=N_e[N_2]$ which implies by Lemma~\ref{lem:exmatch} that $G$ is not well-edge-dominated, which is a contradiction. 

Next, suppose that $G^{xy}_1 \in\{C_7, C_7^*\}$ and label the vertices of $G^{xy}_1$ as $\{v_1, v_2, \dots,  v_7\}$ where $v_1 = z$ and $v_iv_{i+1}\in E(G)$ for $i\in[7]$.  Let $N_1=\{xz, v_2v_3, v_4v_5, v_6v_7\}$ and $N_2=\{xz, v_3v_4, v_5v_6\}.$ Note that $G$ may contain edges between $y$ and vertices in $\{v_4, v_5, v_6\}$, however of each of these vertices is saturated by $N_1$ and $N_2$. Using similar arguments as above, it follows that $N_e[N_1]=N_e[N_2]$ contradicting the fact that $G$ is well-edge-dominated. Therefore, this case cannot occur.

 It follows that for the remainder of the proof we shall assume that $G^{xy}_1$ is bipartite. Furthermore, as our choice of $z$ was arbitrary, we may assume that $G_1^{xz}$ and $G_1^{yz}$ are bipartite. Now relabeling if necessary, we let $z$ be a vertex of $xyz$  such that $G^{xy}_1$  is nontrivial.  We write $V(G^{xy}_1) = A_z \cup B_z$ where $|A_z| \le |B_z|$. We shall assume that $k_{xy} = k$ and write the components of $G^{xy}$ as $G_1^{xy}, \dots, G_k^{xy}$.  Again relabeling $x$ and $y$ if necessary, we may assume that $x$ is adjacent to some vertex of each of $G_2^{xy}, \dots, G_j^{xy}$ where $j \le k$ (reindexing if necessary) and $x$ is not adjacent to any vertex in $G_{j+1}^{xy}\cup \cdots \cup G_k^{xy}$. Thus, $y$ is adjacent to some vertex in each of $G_{j+1}^{xy}, \dots, G_k^{xy}$. Additionally, we may assume that $k>1$ for otherwise the first part of (a) is true and all that remains to show is that if $\deg_G(x) = 2$ or $\deg_G(y) = 2$, then $z$ is detachable in $G_1^{xy}$.
 
Note that $G_1^{yz}$ is bipartite and contains all vertices of $G_2^{xy}\cup \cdots \cup G_j^{xy}$, as well as possibly some vertices from $G^{xy}_1- \{z\}$. Thus,  the graph induced by $\bigcup_{i=2}^j V(G_i^{xy}) \cup \{x\}$, call it $H_x$, is connected and bipartite and we can write $V(H_x) = C_x \cup D_x$ where $|C_x| \le |D_x|$. Similarly, if $j\ne k$, then the graph induced by $\bigcup_{i=j+1}^k V(G_i^{xy}) \cup \{y\}$, call it $H_y$, is connected and bipartite and we can write $V(H_y) = C_y \cup D_y$ where $|C_y| \le |D_y|$. Note that $V(G^{xy}_1) \cup V(H_x) \cup V(H_y)$ is a partition of $V(G)$. 
\vskip5mm
\noindent\textbf{Claim 1:}  $z \in B_z$ 
\vskip2mm
\noindent\emph{Proof:} Assume to the contrary that  $z \in A_z$. Since $G^{yz}_1$ is bipartite, we write $V(G^{yz}_1) = A_x \cup B_x$ where $|A_x| \le |B_x|$. 

Suppose first that $x \in A_x$. First, choose a maximal matching $F_z$ from $G^{xy}_1$ that saturates $A_z$. If $j=2$ and $x$ is not adjacent to any vertex in $G_1^{xy}$, then $G_1^{yz} = H_x$, $A_x = C_x$ and we can find  a maximal matching $F_x$ that saturates $C_x$. On the other hand, if $j>2$ or $x$ is adjacent to a vertex in $G_1^{xy}$, then $x$ is a cut-vertex in $G^{yz}_1$. Further, $H_x$ is a proper connected subgraph of $G_1^{yz}$. Therefore, by Lemma~\ref{lem:subbip}, if $x \in C_x$, we can choose a maximal matching in $H_x$ that saturates $C_x$, and if $x \in D_x$, we can choose a maximal matching in $H_x$ that saturates $D_x$. In either case, we refer to the maximal matching which saturates the partite set of $H_x$ containing $x$ as $F_x$. 

Suppose first that $j<k$, implying that $H_y$ is nontrivial. Again invoking Lemma~\ref{lem:subbip}, $H_y$ is a bipartite subgraph of the well-edge-dominated bipartite graph $G_1^{xz}$ where $y$ is a cut-vertex of $G_1^{xz}$ or $k=j+1$ and $y$ is not adjacent to any vertex in $G_1^{xy}$ or $H_x$. In either case, we can choose a maximal matching in $H_y$ that saturates $y$, call it $F_y$, by Observation~\ref{obs:matchings}. Note that $F_x \cup F_y \cup F_z$ is a maximal matching in $G$ as $F_z$ edge dominates $G_1^{xy}$, $F_x$ edge dominates $H_x$, $F_y$ edge dominates $H_y$ and each vertex of $x, y$, and $z$ is saturated. Moreover, $|F_z \cup F_x \cup F_y| = |A_z| + |C_x| + |C_y|$. On the other hand, if we let $f_z$ be the edge in $F_z$ incident to $z$, and we let $f_x$ be the edge in $F_x$ incident to $x$, then we claim that \[F' = (F_z - \{f_z\}) \cup (F_x - \{f_x\}) \cup F_y \cup \{xz\}\] is also a maximal matching in $G$. To see this, note that $F_z - \{f_z\}$ edge dominates all edges of $G_1^{xy}$ except those edges incident to $z$ as all edges in $G_1^{xy}$ are between a vertex in $A_z$ and a vertex in $B_z$ and $F_z - \{f_z\}$ saturates $A_z - \{z\}$. Similarly, if $x \in C_x$ (resp. $x \in D_x$), then $F_x - \{f_x\}$ edge dominates all edges of $H_x$ except those edges incident to $x$ as all edges in $H_x$ are between a vertex in $C_x$ and a vertex in $D_x$ and $F_x - \{f_x\}$ saturates $C_x - \{x\}$ (resp. $D_x - \{x\}$). However, $|F'| = |F_z \cup F_x \cup F_y| -1$, contradicting the fact that $G$ is well-edge-dominated. 

Therefore, we shall assume $j=k$ and $H_y$ contains only the vertex $y$. If $F_x \cup F_z$ is a maximal matching in $G$, then \[F' = (F_z - \{f_z\}) \cup (F_x - \{f_x\})  \cup \{xz\}\] is also a maximal matching in $G$ for the same reasoning as above, implying that $G$ is not well-edge-dominated. Thus, we shall assume that $F_x \cup F_z$ is not a maximal matching. It follows that the only edges which are not edge dominated by $F_x \cup F_z$ are incident to $y$. In this case, let $e$ be any edge incident to $y$ other than $yz$ or $xy$. Thus, $F_x \cup F_z \cup \{e\}$ is a maximal matching as well as \[F' = (F_z - \{f_z\}) \cup (F_x - \{f_x\})  \cup \{xz,e\}\] which is another contradiction. Hence, this case cannot occur. 

Next, suppose $x \in B_x$. Note that we can assume $|B_x| > |A_x|$ for otherwise $|A_x| = |B_x|$ and we could interchange the roles of $A_x$ and $B_x$ and use all the previous arguments in the case that $x \in A_x$. First, choose a maximal matching $F_z$ from $G_1^{xy}$ that saturates $z$. 

Suppose first that $j =k$ and $H_y$ contains only the vertex $y$. From Lemma~\ref{lem:subbip}, we can choose a maximal matching $M_1$ from $H_x$ that doesn't saturate $x$. Alternatively, we can choose a maximal matching $M_2$ from $H_x$ that does saturate $x$. Using similar arguments as above, one can verify that both $F_z \cup M_1 \cup \{xy\}$ and $(F_z - \{f_z\}) \cup M_2 \cup \{zy\}$ are maximal matchings in $G$ as each of $G_1^{xy}$ and $H_x$ are edge dominated and both sets saturate $x, y$, and $z$. However, $|M_1| = |M_2| = |C_x|$ meaning that we have two different maximal matchings of different cardinality. 

Thus, we shall assume that $j<k$. Note from above that we can assume $y \in B_y$ for otherwise we can interchange the roles of $x$ and $y$ in our above arguments. Moreover, we may assume $|B_y| > |A_y|$ for otherwise $|A_y| = |B_y|$ and we could interchange the roles of $x$ and $y$ as well as interchange the roles of $A_y$ and $B_y$ and use all the previous arguments. By Lemma~\ref{lem:subbip}, we can choose a maximal matching $M_1$ from $H_x$ that doesn't saturate $x$ as well as a maximal matching $N_1$ from $H_y$ that doesn't saturate $y$. Similarly, we can choose a maximal matching $M_2$ from $H_x$ that does saturate $x$ as well as a maximal matching $N_2$ from $H_y$ that does saturate $y$. It follows that $F_z \cup M_1 \cup N_1 \cup \{xy\}$ and $F_z \cup M_2 \cup N_2$ are maximal matchings in $G$ with different cardinality, another contradiction. Hence, this case cannot occur. 
\smallqed
\vskip2mm

It follows that for the remainder of the proof that we may assume $z \in B_z$. Moreover, as our choice of $z$ from $\{x, y, z\}$ was arbitrary, we may assume that if $H_x$ is nontrivial, then $x \in B_x$ (in $G_1^{yz}$) and if $H_y$ is nontrivial, then $y \in B_y$ (in $G_1^{xz}$). Furthermore, from previous arguments, we may assume $|B_z| > |A_z|, |B_x| > |A_x|$ and if $H_y$ is nontrivial, then $|B_y| > |A_y|$.

Suppose that $G_2^{xy}$ exists. First, we argue that $H_y$ is trivial. Suppose to the contrary that $H_y$ is not trivial. Choose a maximal matching $F_z$ from $G_1^{xy}$ that saturates $z$, a maximal matching $F_x$ from $H_x$ that saturates $x$, and a maximal matching $F_y$ from $H_y$ that saturates $y$. Alternatively, choose a maximal matching $M_1$ from $H_x$ that doesn't saturate $x$ and a maximal matching $N_1$ from $H_y$ that doesn't saturate $y$. One can easily verify that  $F_z \cup F_x \cup F_y$ and $F_z \cup M_1 \cup N_1 \cup \{xy\}$ are maximal matchings of different cardinality as $|F_x| = |M_1|$ and $|F_y| = |N_1|$. Therefore, we may assume that $H_y$ is trivial. Next, we show that $\deg_G(y) \ge 3$. Assume to the contrary that $\deg_G(y) = 2$. Then we can choose a maximal matching $M_1$ in $H_x$ that saturates $x$ and a maximal matching  $M_2$ in $G_1^{xy}$ that saturates $z$ and $M_1 \cup M_2$ is a maximal matching in $G$. On the other hand, we can choose a maximal matching $M_3$ in $H_x$ that doesn't saturate $x$ and a maximal matching $M_4$ in $G_1^{xy}$ that doesn't saturate $z$ and $M_3 \cup M_4 \cup \{xz\}$ is a maximal matching in $G$. 

Finally, we show that if either $x$ or $y$ have degree $2$ in $G$, or $G_2^{xy}$ exists, then $z$ is detachable in $G_1^{xy}$. Suppose to the contrary that $G_1^{xy} - z$ is not well-edge-dominated. Thus, we can find two minimal edge dominating sets $F_1$ and $F_2$ in $G_1^{xy}-z$ where $|F_1| < |F_2|$. Note that if $\deg_G(x) = 2$ (or $\deg_G(y)=2$), then both $F_1 \cup \{yz\}$ and $F_2 \cup \{yz\}$ (or $F_1 \cup \{xz\}$ and $F_2 \cup \{xz\}$) are minimal edge dominating sets in $G$, which is a contradiction. Therefore, we shall assume that $\deg_G(x)\ge 3$ and $\deg_G(y)\ge 3$ and $G_2^{xy}$ exists. By assumption $H_y$ is trivial. Moreover, $H_x'$ obtained from $H_x$ by removing the vertices in $G_1^{xy}$ is a well-edge-dominated graph with cut-vertex $x$. If $H_x \ne H_x'$, then we can write $H_x - x = G_2^{xy} \cup \cdots \cup G_j^{xy} \cup J_1 \cup \cdots J_m$ where each $J_i$ is a subgraph of $G_1^{xy}$. One can easily verify that $H_x'$ is well-edge-dominated as $H_x$ is assumed to be well-edge-dominated and we know $H_x' = H_x - N_e[I]$ where $I$ is a maximum matching in  $J_1 \cup \cdots \cup J_m$ is therefore well-edge-dominated. Thus, by Lemma~\ref{lem:subbip}, there exists a maximal matching $M$ in $H_x'$ which saturates $x$. It follows that $F_1 \cup M \cup \{yz\}$ and $F_2 \cup M \cup \{yz\}$ are minimal edge dominating sets in $G$, another contradiction. 

\end{proof}

Now we proceed by considering  well-edge-dominated graphs $G$ containing exactly one triangle $xyz$ and the possibilities of the degree sequence of $x, y$, and $z$ as well as the number of nontrivial components in the graph induced by $V(G) - \{x, y, z\}$. Throughout the remainder of the section, we define $G'$ to be precisely this graph induced by $V(G) - \{x, y, z\}$. We first consider when $G'$ contains a nontrivial component, $G_i'$, such that for some vertex $a \in \{x, y, z\}$, each vertex $v$ of $G_i'$ is adjacent to $a$ or $v$ is not adjacent to any vertex of $\{x, y, z\}$. We use the following lemma shown in \cite{AKR-2022}.

\begin{lemma}\cite{AKR-2022}\label{lem:Krs} If $2 \le r < s$, then $K_{r, s}$ is not well-edge-dominated. 
\end{lemma}

\begin{theorem}\label{thm:houseonfire} Let $G$ be a well-edge-dominated graph containing exactly one triangle $xyz$ such that the graph induced by $V(G) - \{x, y, z\}$ contains a nontrivial component with no vertex adjacent to either $x$ or $y$. Then $G \in \mathcal{F}\cup \mathcal{T}$. 
\end{theorem}

\begin{proof} 
Suppose a counterexample exists and let $G$ be such a counterexample of smallest order. Consider the graph $G'$ constructed from $G$ by removing $\{x, y, z\}$. Note that $G'$ may or may not be disconnected. Suppose first $G'$ contains only one nontrivial component, say $G_1'$.  By assumption $G_1'$ contains no vertex adjacent to either $x$ or $y$. If $G'$ contains trivial components, then each trivial component of $G'$ is a leaf of $G$, adjacent to one of $x, y$, or $z$. Note that $z$ cannot be a support vertex of $G$ as this would imply that the component of $G-N_e[xy]$ that contains $z$, call it $G_z$, is a well-edge-dominated bipartite graph with $V(G_z) = A_z \cup B_z$ and $|A_z|\le |B_z|$ where $z \in A_z$, contradicting Theorem~\ref{thm:TRI1}. So if $G'$ contains trivial components, we may assume that either $x$ or $y$ is a support vertex in $G$, but not both according to Theorem~\ref{thm:TRI1}. However, if either $x$ or $y$ is a support vertex, then neither (a) nor (b) in Theorem~\ref{thm:TRI1} is true, which is a contradiction. It follows that $G' = G_1'$ and $G \in \mathcal{T}$. 

Therefore, for the remainder of the proof, we may assume that $G'$ contains at least two nontrivial components. Write $G' = G_1'\cup G_2'\cup \cdots \cup G_k'$ where we may assume $z$ is adjacent to a vertex in each of $G_1', \dots, G_{\alpha}'$ and $z$ is not adjacent to any vertex in $G_{\alpha+1}', \dots, G_k'$. Further, if we let $G_z$ be the component of $G-N_e[xy]$ that contains $z$, then $V(G_z) = V(G_1')\cup \cdots\cup V(G_{\alpha}') \cup \{z\}$. If need be relabel $G_1', \dots, G_{\alpha}'$ such that for some $r \in [\alpha]$, each component of $G_i'$ has no vertex adjacent to either $x$ or $y$ for $i \in [r]$ and $G_i'$ has at least one vertex adjacent to either $x$ or $y$ for $r+1 \le i \le \alpha$. By assumption, we may assume $G_1'$ has no vertex adjacent to either $x$ or $y$. 

 Note that we know $G_1'$ is a well-edge-dominated bipartite graph by Theorem~\ref{thm:TRI1}. Moreover, we know that  $G_z$ is a well-edge-dominated bipartite graph where the bipartition is $A_z\cup B_z$ with $|A_z| < |B_z|$ and $z \in B_z$ by Theorem~\ref{thm:TRI1}. It follows that $G_1'$ is not trivial for otherwise $z$ is a support vertex in $G_z$ contradicting Theorem~\ref{thm:TRI1}. Let $w \in A_z\cap V(G_1')$ be a neighbor of $z$ in $G$. Choose a  maximal matching $F_1$ in $G_1'$ which saturates $w$. Now consider the graph $H= G-N_e[wz]$ and let $H_{xy}$ be the component of $H$ containing $x$ and $y$. We know each component of $H$ is well-edge-dominated. 
 
 We first argue that $H_{xy}$ is not in $\{C_5, C_7, C_7^*\}$. Suppose first that $H_{xy} = C_5$ and write the $5$-cycle as $xyabcx$. It follows that $abc$ is a path in some $G_i'$. Moreover, $G_i' = P_3$ for otherwise $H_{xy}$ does not have order $5$. However, $J = G - N_e[ab]$ is a well-edge-dominated graph containing the triangle $xyz$ where $G_1'$ is still a component of $J$ and $V(J) - \{x, y, z\}$ contains a nontrivial component with no vertex adjacent to either $x$ or $y$. However, $J \not\in \mathcal{T} \cup \mathcal{F}$ as $\deg_J(x) =3$, $\deg_G(y) = 2$,  and $\deg_J(z) \ge 3$, contradicting the minimality of $n(G)$. Next, suppose $H_{xy} \in \{C_7, C_7^*\}$ and write the $7$-cycle as $xyabcdex$. As above, this implies that some $G_i'$ is the path $abcde$ which yields $J = G-N_e[ab]$ is a smaller counterexample. Thus, we may assume that $H_{xy} \not\in \{C_5, C_7, C_7^*\}$ and $H_{xy}$ is bipartite. 
 
 Suppose first that $H_{xy} =K_2$. In this case, $\deg_G(x) = \deg_G(y) = 2$ and by Theorem~\ref{thm:TRI1} (a), $z$ is detachable in $G_z $ and $G\in \mathcal{T}$. Therefore, we shall assume that $H_{xy} \ne K_2$. Suppose first that there exists a maximal matching $M_1$ in $H_{xy}$ that saturates $x$ and another maximal matching $M_2$ in $H_{xy}$ that does not saturate $x$. As $H_{xy}$ is well-edge-dominated, $|M_1| = |M_2|$. It follows that when we partition $H_{xy}$ as $A_{xy}\cup B_{xy}$ where $|A_{xy}| \le |B_{xy}|$ that $x \in B_{xy}$ and $y \in A_{xy}$. Thus, $y$ is saturated in both $M_1$ and $M_2$. Let $f_y$ be the edge in $M_1$ incident to $y$. Pick a maximal matching $F_i$ in $G_i'$ for $1 \le i \le r$. If $f_y \ne xy$, then  $(M_1 - \{f_y\}) \cup  \bigcup_{i=1}^r F_i \cup \{yz\}$ is a maximal matching in $G$ as all vertices of the triangle are saturated and all edges of $G'$ are dominated. On the other hand, $M_2 \cup  \bigcup_{i=2}^r F_i \cup \{xz\}$ is also a maximal matching, which is a contradiction. Therefore, we shall assume that the only maximal matching $M_1$ in $H_{xy}$ that saturates $x$ must contain $xy$, i.e. $x$ is a leaf in $H_{xy}$ and $\deg_G(x) = 2$. In this case, $(M_1 - \{xy\}) \cup  \bigcup_{i=1}^r F_i \cup \{yz\}$ and $M_2 \cup  \bigcup_{i=1}^r F_i \cup \{xz\}$ are maximal matchings, yet another contradiction.  Similarly, if there exist two maximal matchings in $H_{xy}$, one that saturates $y$ and the other that doesn't, we reach another contradiction. Thus, we may assume that any maximal matching in $H_{xy}$ saturates both $A_{xy}$ and $B_{xy}$, implying $|A_{xy}| = |B_{xy}|$. As the only well-edge-dominated balanced bipartite graphs are complete bipartite graphs, it follows that $H_{xy} = K_{n,n}$. Moreover, we may assume $n >1$ as $H_{xy} \ne K_2$. This implies that $z$ is not adjacent to any vertex in $H_{xy}$ (for otherwise $G$ contains two triangles) and $G'$ consists of $G_1', \dots, G_r'$ together with the component $K_{n-1, n-1}$. Now when we consider $G - N_e[xz]$, the component containing $y$ is isomorphic to $K_{n, n-1}$, which is not well-edge-dominated unless $n=2$ by Lemma~\ref{lem:Krs}. It follows that $H_{xy}= K_{2, 2}$ and $r=\alpha$.

Write $H_{xy} = xyst$. We need to show that $z$ is strongly detachable in $G_z$. It is clear that $G$ has a maximal matching of size $|A_z| + 2$. Also note that $G-N_e[xy]$ has two components; namely, $G_z$ and $st$. Since $G_z$ is well-edge-dominated, $z \in B_z$, and $z$ is a cut-vertex in $G_z$, $G_z - z$ is well-edge-dominated by Lemma~\ref{lem:subbip2}. All that remains is to show that each vertex in $N_{G_z}(z)$ is a support vertex in $G_z - z$. Suppose to the contrary that $w \in G_z - z$ is not a support in $G_z - z$. We claim that we can find a minimal edge dominating set of $G_z-z$ which does not saturate $w$. Indeed, let $U_w = N_{G_z}(w) - \{z\}$ and $V_w = N_{G_z}(U_w) - \{w\}$. For each vertex $u \in U_w$, we can pick an edge $uv$ such that $v \in V_w$, and call the resulting set $F$. Now let $G^0$ be the bipartite subgraph of $G_z - z$ containing edges not dominated by $F$ ($G_0$ is possibly empty). Let $M$ be a maximal matching for $G_0$. We claim that $F \cup M$ is a minimal edge dominating set of $G_z - z$. To see this, note that each edge in $M$ is its own private neighbor with respect to $F \cup M$ and each edge $uv \in F$ has its own private neighbor of the form $aw$. Since $G_z - z$ is well-edge-dominated, $|F\cup M| = |A_z|$. Now $F' = F \cup M \cup \{xz, xy, xt\}$ is a minimal edge dominating set of $G$ as the private neighbor of $xz$ is $wz$, the private neighbor of $xy$ is $ys$, the private neighbor of $xt$ is $st$ and every edge in $F\cup M$ still has a private neighbor in $G$ (as no edge incident to $z$ in $G_z$ is in $F'$). However, $|F'| = |A_z| + 3$ which is a contradiction. Thus, each vertex in $N_{G_z}(z)$ is a support vertex  in $G_z-z$ meaning $z$ is strongly detachable, so $G \in \mathcal{F}$.
\end{proof}

Based on the previous result and Lemma~\ref{lem:noleaves}, we may assume that $G'$ contains a nontrivial component and every nontrivial component of $G'$ contains a pair of vertices $u$ and $v$ such that $u$ is adjacent to some $a \in \{x, y, z\}$, and $v$ is adjacent to some vertex in $\{x, y, z\} - \{a\}$. Moreover, one of the following must be true:
\begin{itemize}
\item There exists some edge $ab$ on $xyz$ such that $G- N_e[ab]$ contains exactly one nontrivial component, or
\item For any edge $f$ on $xyz$, $G - N_e[f]$ does not contain exactly one nontrivial component. 
\end{itemize}
Furthermore, when we consider the possible degree sequences of $y, x$, and $z$ in $G$ (where we may relabel so that $\deg_G(y) \le \deg_G(x) \le \deg_G(z)$), we have one of three scenarios: (1) the degree sequence is $(2, 2, k)$ for some $k \ge 3$, or (2) the degree sequence is $(2, k, \ell)$ for some $3 \le k \le \ell$, or (3) $\deg_G(y) \ge 3$. Note that in the first case, $G - N_e[xy]$ contains exactly one nontrivial component, call it $G_0$, as it cannot contain only isolates. However, we have assumed that some vertex of $G_0$ is adjacent to either $x$ or $y$ so this case cannot occur. Moreover, if every vertex of $x, y$, and $z$ has degree at least $3$, and for some edge on $xyz$, say $xy$, $G - N_e[xy]$ contains only trivial components, then $\deg_G(z) = 2$, which is a contradiction. Hence, there are only three cases to consider:
\begin{enumerate}
\item[(1)] Exactly one vertex of $\{x, y, z\}$ has degree $2$.
\item[(2)] For any $f$ on $xyz$, $G- N_e[f]$ contains at least two nontrivial components and each vertex on $xyz$ has degree at least $3$ in $G$. 
\item[(3)] There exists some edge $f$ on $xyz$ such that $G - N_e[f]$ contains exactly one nontrivial component and each vertex on $xyz$ has degree at least $3$ in $G$. 
\end{enumerate}
The next three results address each of the above (in the order given).

\begin{theorem}\label{thm:degree2} Let $G$ be a connected, well-edge-dominated graph containing exactly one triangle $xyz$. If $\deg_G(y)=2$, $\deg_G(x) \ge 3$ and $\deg_G(z) \ge 3$, then $G \in\{ \mathcal{H}, \mathcal{DH}\}$. 
\end{theorem}

\begin{proof} Let $G_1', \dots, G_k'$ be the components of $G'$. We may assume that each $G_i'$ contains a vertex adjacent to either $x$ or $z$, or both, as $\deg_G(y) = 2$. Suppose first that $G_1'$ contains a vertex adjacent to $z$, but no vertex adjacent to $x$. By Theorem~\ref{thm:TRI1}, we know that the component of $G-N_e[xy]$ containing $z$ is bipartite with bipartition $A_z\cup B_z$ where $|A_z|<|B_z|$ and $z \in B_z$. Thus, $G_1'$ is not trivial. Furthermore, by Theorem~\ref{thm:houseonfire}, we know $G \in \mathcal{F}\cup \mathcal{T}$. Since $\deg_G(y)=2$, it must be that $\deg_G(x) = 2$, which is a contradiction. Therefore, no such component of $G'$ contains a vertex adjacent to $z$ but no vertex adjacent to $x$. A similar argument can be used to show that no component of $G'$ contains a vertex adjacent to $x$ but no vertex adjacent to $z$. So we shall assume that each component in $G'$ contains a vertex adjacent to $x$ and a vertex adjacent to $z$. It follows that every component $G_i'$ is nontrivial for otherwise $G$ contains two triangles. Now let $H = G-N_e[xz]$ which is bipartite (not necessarily connected) and well-edge-dominated. Therefore, we write $V(H) = A \cup B$ where $|A| \le |B|$. 

Assume first that $|A|<|B|$. We claim that $z$ does not have a neighbor $w \in B$. Indeed, suppose this is not the case and $w \in B$ is adjacent to $z$. It follows that we can choose two different maximal matchings in $H$ of size $|A|$, $F_1$ that saturates $w$ and $F_2$ that does not, by Observation~\ref{obs:matchings}.  However,  $F_1 \cup \{xz\}$ and $F_2 \cup \{wz, xy\}$ are two different maximal matchings in $G$, contradicting that $G$ is equimatchable. So we may assume that the only neighbors of $z$ in $H$ are in $A$. Similarly, the only neighbors of $x$ in $H$ are in $A$. Letting $F$ be any maximal matching in $H$ (which necessarily saturates $A$ as $H$ is well-edge-dominated), assume $f_w \in F$ is incident to $w$, a neighbor of $z$, and assume $f_t \in F$ is incident to $t$, a neighbor of $x$. Note that $f_w \ne f_t$ as $G$ contains only one triangle. Now $F \cup \{xz\}$ and $(F - \{f_w, f_t\}) \cup \{zw, xt\}$ are two maximal matchings in $G$, another contradiction. Therefore, it cannot be that $|A|<|B|$ and we may assume that $|A|=|B|$. 

It follows that every component of $H$ is a complete bipartite graph. Suppose one component of $H$, call it $H_1$, is isomorphic to $K_{n,n}$ where $n \ge 3$. Without loss of generality, we may assume that some $w\in V(H_1)$ is a neighbor of $z$ (as some vertex of $H_1$ is adjacent to $x$ or $z$). Then $G-N_e[\{zw, xy\}]$ contains the component $K_{n, n-1}$  which is not well-edge-dominated by Lemma~\ref{lem:Krs}. Thus, every component of $H$ is either $K_2$ or $C_4$. Moreover, we have assumed that each $K_2$ or $C_4$ component contains a vertex adjacent to $x$ and a vertex adjacent to $z$. Suppose first that two components of $H$ are isomorphic to $K_2$, say  $w_1t_1$ and $w_2t_2$ where $w_i$ is adjacent to $z$ and $t_i$ is adjacent to $x$ for $i \in [2]$. For any maximal matching $M$ in $H$, $M \cup \{xz\}$ and $(M- \{w_1t_1, w_2 t_2\}) \cup \{w_1z, xt_2\}$ are two maximal matchings in $G$, which is a contradiction. Therefore, we may assume that at most one component of $H$ is a $K_2$. Suppose that $H_1 = C_4 = wrst$ and $H_2 = C_4 = mn\ell o$ are two components in $H$ where $z$ is adjacent to $w$ and $x$ is adjacent to $m$. Pick a maximal matching $M$ from $H$ and let $f_w \in M$ be incident to $w$ and $f_m \in M$ be incident to $m$. Then $M \cup \{xz\}$ and $(M - \{f_w, f_m\}) \cup \{zw, xm\}$ are maximal matchings in $G$. Furthermore, a similar argument can be used to show that $H$ cannot contain two components, one isomorphic to $C_4$ and another isomorphic to $K_2$. It follows that $H \in \{C_4, K_2\}$. If $H = K_2$, then $G$ is the house graph. So assume that $H= C_4 = wrst$ where $w$ is adjacent to $z$ and at least one of $\{r, s, t\}$ is adjacent to $x$.

Assume first that $x$ is adjacent to $r$ or $t$. With no loss of generality, we may assume that $x$ is adjacent to $r$. If $\deg_G(s) = \deg_G(t)=2$, then $G$ is the dream house graph. So assume $G$ contains more edges. Since $\deg_G(y)=2$ and $\{rz, wx\} \cap E(G) = \emptyset$, the only additional edges in $G$ are incident to either $s$ or $t$. If $zs \in E(G)$, then $G$ contains a mazimal matching of size $3$ yet $\{xy, wr, wt, wz\}$ is a minimal edge dominating set as $zs$ is the private edge neighbor of $wz$, $sr$ is the private edge neighbor of $wr$, $st$ is the private edge neighbor of $wt$, and $xy$ is its own private edge neighbor. This cannot be as $G$ is assumed to be well-edge-dominated. On the other hand, $s$ is not adjacent to $x$ either as $G$ contains exactly one triangle. Hence, $\deg_G(s) =2$ and similar arguments show that $\deg_G(t) = 2$, another contradiction. Thus, we may assume that $x$ is not adjacent to $r$, and (by similar logic) $x$ is not adjacent to $t$. It follows that $x$ is adjacent to $s$ as some vertex of $\{r, s, t\}$ is adjacent to $x$. Hence, $\deg_G(r) = \deg_G(t) = 2$ as $G$ contains only one triangle. In this case, $G$ contains a matching of size $3$ and yet $\{zw, xs\}$ is an edge dominating set, another contradiction. Having exhausted all possibilities, the result follows. 
\end{proof}

\begin{theorem}\label{thm:houseonfire2} Let $G$ be a connected, well-edge-dominated graph containing exactly one triangle $xyz$. If $G-N_e[f]$ contains at least two nontrivial components for any edge $f$ on $xyz$, and each of $x, y$, and $z$ have degree at least $3$ in $G$, then $G \in \mathcal{F}$. 
\end{theorem}

\begin{proof} 
Consider the graph $G'$ constructed from $G$ by removing $\{x, y, z\}$. Note that $G'$ is disconnected by assumption and contains at least two nontrivial components. Write $G' = G_1'\cup G_2'\cup \cdots \cup G_k'$ where we may assume $z$ is adjacent to a vertex in each of $G_1', \dots, G_{\alpha}'$ and $z$ is not adjacent to any vertex in $G_{\alpha+1}', \dots, G_k'$. Further, if we let $G_z$ be the component of $G-N_e[xy]$ that contains $z$, then $V(G_z) = V(G_1')\cup \cdots\cup V(G_{\alpha}') \cup \{z\}$. By assumption, $G-N_e[xy]$ contains at least two nontrivial components so $G_{\alpha+1}'$ also exists. From Theorem~\ref{thm:houseonfire}, we may assume that each vertex in $\{x, y\}$ is adjacent to a vertex in each of $G_{\alpha+1}', \dots, G_k'$. Moreover, we may assume by Theorem~\ref{thm:houseonfire} that each component in  $G_1', \dots, G_{\alpha}'$  has at least one vertex adjacent to either $x$ or $y$.

Now, we have assumed that each component of $G_1', \dots, G_{\alpha}'$ contains a vertex adjacent to $z$ and contains a vertex adjacent to either $x$ or $y$, or both.  By simply interchanging the roles of $x, y$, and $z$, we may also assume that some component in $G_1', \dots, G_{\alpha}'$ contains no vertex adjacent to $y$ and some component in $G_1', \dots, G_{\alpha}'$ contains no vertex adjacent to $x$. So we relabel the components $G_1', \dots, G_{\alpha}'$ as \[H_{xy}^1, \dots, H_{xy}^j, H_{xz}^1, \dots, H_{xz}^k, H_{yz}^1, \dots, H_{yz}^{\ell}, H_{xyz}^1, \dots, H_{xyz}^m\] such that for component $H_I^r$, for each vertex $u \in I$, there exists a vertex in $H_I^r$ that is adjacent to $u$.

Now let $f_y$ be an edge incident to $y$ and a vertex $u$ of $H_{yz}^1$. We claim that $J = G- N_e[f_y]$ is bipartite and connected. It is clear that each component of $G'$ has a vertex adjacent to either $x$ or $z$. Thus, $J$ is connected. To see that $J$ is bipartite, note that $H_{xy}^1$ and $H_{xz}^1$ exist and contain at least two vertices while $H_{yz}^1$ contains at least one vertex other than $u$. Hence, $n(J) \ge 7$ and if $J$ is not bipartite, then $J \in \{C_7, C_7^*\}$. Moreover, since $x$ has degree $3$ in $J$, $J = C_7^*$. It follows that $H_{xy}^1$ contains exactly two vertices, say $a$ and $b$, where $a$ is adjacent to $x$ and $b$ is a leaf in $J$, which cannot be. It follows that $J$ is bipartite. Therefore, we may write $V(J) = A_J \cup B_J$ where $|A_J| \le |B_J|$, and  $x$ and $z$ are in different partite sets of $J$. On the other hand, we know that $G_z$ is a well-edge-dominated, and therefore equimatchable, bipartite graph with cut-vertex $z \in B_z$ where $V(G_z) = A_z \cup B_z$ and $|A_z|<|B_z|$. It follows that $H_{xz}^i$ for any $i\in [k]$ is bipartite where $V(H_{xz}^i) = A_{xz}^i \cup B_{xz}^i$ with $|A_{xz}^i| \le |B_{xz}^i|$. Moreover, $H_{xz}^i$ is well-edge-dominated as it is a component of $G- N_e[xz]$. By Lemma~\ref{lem:subbip2},  $z$ is only adjacent to vertices in $A_{xz}^1$.  However, when we consider $G_x$ as the component in $G-N_e[yz]$ that contains $x$ and $H_{xz}^i$ is a component of $G_x-x$, then $x$ can only be adjacent to vertices in $B_{xz}^i$ as $J$ is bipartite. Moreover, when we interchange the roles of $x$ and $z$, we know that $x \in B_x$ and therefore with respect to $G_x-x$, the roles of $A_{xz}^i$ and $B_{xz}^i$ must interchange as $x$ and $z$ are in the bipartite graph $J = G- N_e[f_y]$. This is only possible if $|A_{xz}^i| = |B_{xz}^i|$ and $z$ is only adjacent to vertices in $B_{xz}^i$ and $x$ is only adjacent to vertices in $A_{xz}^i$ (or vice versa). Additionally, this same argument shows that $H_{xyz}^1$ does not exist. To see this, suppose $H_{xyz}^1$ exists. Thus, by the above argument, at least two vertices in $\{x, y, z\}$ would be adjacent to vertices in the same partite set of $H_{xyz}^1$. Relabeling if necessary, we may assume that $y$ and $z$ are both adjacent to vertices in $A_{xyz}^1$ (or $B_{xyz}^1$). If there exists a neighbor of $x$, call it $r$, that is in $H_{xz}^1 \cup H_{xy}^1 \cup H_{xyz}^2$, then $G - N_e[xr]$ contains an odd cycle involving $y$, $z$, and vertices in $H_{xyz}^1$. However, this is only possible if the component of $G - N_e[xr]$ containing $y$ and $z$ is $C_5, C_7$, or $C_7^*$. This would in turn imply that $H_{xz}^1$, $H_{xy}^1$, and $H_{xyz}^2$ do not exist, a clear contradiction. Thus, the only neighbors of $x$ other than $z$ and $y$ are in $H_{xyz}^1$. However, in this case, $G - N_e[xz]$ contains only one nontrivial component which cannot be. Thus, we may assume that $H_{xyz}^1$ does not exist. 

Note that since $H_{xz}^i$ is equimatchable, it follows that $H_{xz}^i \cong K_{n,n}$ for some $n \in \mathbb{N}$. Suppose first  that  $n\ge 3$. We could pick any edge $f_z$ incident to $z$ and a vertex in $H_{xz}^i$ and then $G - N_e[\{f_z, xy\}]$ contains a component which is  isomorphic to a subgraph $K_{n, n-1}$ which is not well-edge-dominated by Lemma~\ref{lem:Krs}. Therefore, this case cannot occur. It follows that $H_{xz}^i \in \{K_2, C_4\}$. Furthermore, this is true for every component of the form $H_{xy}^i$ or $H_{zy}^i$. We find two different maximal matchings in $G$ of different cardinality. First, choose a maximal matching $F$ from $G'$ which saturates every neighbor of $z$ in $G'$ so that $F \cup \{xy\}$ is maximal in $G$. Next, pick  $f_x$ to be an edge incident to $x$ and a vertex in $H_{xy}^1$, pick $f_z$ incident to $z$ and a vertex in $H_{xz}^1$, and pick $f_y$ incident to $y$ and a vertex in $H_{yz}^1$. For each graph in $\{H_{xy}^1, H_{xz}^1, H_{yz}^1\}$ which is isomorphic to $C_4$, pick an additional edge that is not adjacent to $f_x, f_y$, or $f_z$ and call the resulting set $M$. Finally, let $N$ be any maximal matching in $G' -(H_{xy}^1\cup H_{xz}^1\cup H_{yz}^1)$. It follows that $M \cup N \cup \{f_x, f_y, f_z\}$ is also a maximal matching in $G$ of cardinality $|F|$, which is a contradiction. 

\end{proof}

\begin{theorem}\label{thm:crystal} Let $G$ be a connected, well-edge-dominated graph containing exactly one triangle $xyz$. If $G-N_e[xy]$ contains one nontrivial component and each vertex on the triangle has degree at least $3$ in $G$, then $G = Cr$. 
\end{theorem}

\begin{proof} Consider the graph $G'$ constructed from $G$ by removing $\{x, y, z\}$ and let $G_1', \dots, G_k'$ be the components of $G'$. By assumption, there exists some vertex in each of $G_i'$ adjacent to $z$. From Theorem~\ref{thm:TRI1}, we know that $G-N_e[xy]$ is bipartite with bipartition $A_z \cup B_z$ such that $z \in B_z$. It follows that no component of $G'$ is an isolate for this would imply that  $z$ is a support vertex in $G-N_e[xy]$, contradicting the fact that $z \in B_z$. Furthermore, by Theorem~\ref{thm:houseonfire}, we may assume that each $G_i'$ contains a vertex adjacent to $x$ or $y$. By Theorem~\ref{thm:TRI1}, $G-N_e[yz]$ and $G-N_e[xz]$ are bipartite and in the component of $G-N_e[yz]$ containing $x$, $x$ is not a support vertex. Similarly, in the component of $G-N_e[xz]$ containing $y$, $y$ is not a support vertex. It follows  that we can enumerate the components of $G'$ as $G_1', \dots, G_{k}'$ where each of  $G_1', \dots, G_j'$ are bipartite, nontrivial, and contain a vertex adjacent to $x$, each of $G_{j+1}', \dots, G_{k}'$ are bipartite, nontrivial, and contain a vertex adjacent to $y$.

For each of $G_1', \dots, G_{k}'$, we partition the vertex set as $A_i'\cup B_i'$ where $|A_i'| \le |B_i'|$. Since $z \in B_z$ we know by Lemma~\ref{lem:subbip2} that for each $ i  \in[k]$ that $z$ is adjacent to some vertex in $A_i'$. Suppose first that $x$ is adjacent to some vertex $w \in B_i'$ for $i \in [j]$ and $|A_i'| < |B_i'|$. Since $G_i'$ is well-edge-dominated, there exist matchings $M_i$ and $M_i'$ in $G_i'$ where $M_i$ saturates $w$ and $M_i'$ does not saturate $w$. For all other $\alpha \in \{1, 2, \dots, i-1, i+1, \dots, k\}$, choose a maximal matching $M_{\alpha}$ in $G_{\alpha}'$. Note that $\bigcup_{\alpha\ne i} M_{\alpha} \cup M_i' \cup \{wx, yz\}$ is a maximal matching in $G$ as it contains a maximal matching in $G'$ and saturates $x, y$, and $z$. On the other hand, $\bigcup_{\alpha\ne i} M_{\alpha} \cup M_i \cup \{xy\}$ is a maximal matching as every vertex from $\cup_{i=1}^k A_i'$ is saturated as well as $x$ and $y$ are saturated. However, these two maximal matchings have different cardinality, which is a contradiction. Thus, $x$ only has neighbors in $A_i'$ where $|A_i'|<|B_i'|$ or $x$ has neighbors in $G_i'$ where $|A_i'| = |B_i'|$. Using a similar argument, we may assume that $y$ only has neighbors in $A_i'$ where $|A_i'|<|B_i'|$ or $y$ has neighbors in $G_i'$ where $|A_i'| = |B_i'|$. 

Now choose a maximal matching $M_i$ for $i \in [k]$. Thus, each $M_i$ saturates $A'_i$ for $i \in [k]$. Suppose first that we can find $3$ edges $\{a_1b_1, a_2b_2, a_3b_3\} \subset \bigcup_{i=1}^{k} M_i$ such that $a_1$ is adjacent to $z$, $a_2 \in A_i'$ for some $i \in [j]$ and $a_2$ is adjacent to $x$, and $a_3 \in A_i'$ for some $j+1 \le i \le k$ and $a_3$ is adjacent to $y$. Then $\bigcup_{i=1}^{k}M_i \cup \{xy\}$ and $\left(\bigcup_{i=1}^{k} M_i - \{a_1b_1, a_2b_2, a_3b_3\}\right) \cup \{za_1, xa_2, ya_3\}$ are maximal matchings in $G$. Thus, we may assume that every neighbor of $x$ in $G_i'$ is in $B_i'$ for $i \in [j]$ or every neighbor of $y$ in $G_i'$ is in $B_i'$ for $j+1\le i \le k$. Without loss of generality, assume that every neighbor of $x$ in $G_i'$ is in $B_i'$ for $i \in [j]$. It follows that $G_i' = K_{n,n}$ for some $n \ge 1$ and $i \in [j]$. Now suppose that $G_1' = K_{n, n}$ where $n \ge 3$. If we assume that $t \in V(G_1')$ is a neighbor of $z$, then $G_1' - t$ is a component of $G- N_e[\{xy, tz\}]$, which is a contradiction as $K_{n, n-1}$ is not well-edge-dominated when $n \ge 3$. Thus, $G_i' \in \{K_2, C_4\}$ for each $i \in [j]$. 

Next, suppose $y$ is adjacent to some vertex in $A_{j+1}'$, let $w$ be a neighbor of $x$ in $G_1'$, and let $t$ be a neighbor of $z$ in $G_1'$. It follows that the component of $G-N_e[wx]$ containing $z$, call it $H$, is well-edge-dominated with girth at least $4$ and contains an odd cycle involving $z, y$, and vertices from $G_{j+1}'$. Thus, $H \in \{C_5, C_7, C_7^*\}$, but this cannot be as $t$ is a neighbor of $z$ in $H$ which is not on the odd cycle containing $z, y$, and vertices from $G_{j+1}'$. Therefore, this case cannot occur and we assume that either $G_1'$ is the only component in $G'$, or every neighbor of $y$ in $G_i'$ is in $B_i'$ for $j+1 \le i \le k$.

Suppose first that every neighbor of $y$ in $G_i'$ is in $B_i'$ for $j+1 \le i \le k$ and $k \ne 1$. As above, this implies that $G_i' \in \{K_2, C_4\}$ for each $j+1 \le i \le k$. Pick a maximal matching $M_i$ from each $G_i'$ for $i \in [k]$, which we know saturates all vertices from $G_i'$. Let $Y\subset \bigcup_{i=1}^{k}M_i$ be those edges that are incident to a neighbor of $y$. Since $G$ contains only one triangle, this implies that no edge in $Y$ is incident to $x$. Now consider $J = G-N_e[Y]$ which contains the triangle $xyz$ where $\deg_J(y) = 2$, $\deg_J(x) \ge 3$ and $\deg_J(z) \ge 3$. It follows that $J \in \{\mathcal{H}, \mathcal{DH}\}$ which means $y$ is adjacent to every vertex of $B_i'$ for $2 \le i \le k$ for otherwise $z$ would have degree $4$ or more in $J$. On the other hand, if we let $X \subset \bigcup_{i=1}^k M_i$ be those edges that are incident to a neighbor of $x$, then by the same argument $G - N_e[X]$ is either $\mathcal{H}$ or $\mathcal{DH}$ so $x$ is adjacent to every vertex of $B_i'$ for $1 \le i \le k-1$. But $x$ and $y$ can have no common neighbor other than $z$ so this case cannot occur unless $j=1$ and $k=2$. Moreover,  $G_1'=K_2=G_2'$ for otherwise $y$ has degree $4$ in $J$ or $x$ has degree $4$ in $G- N_e[X]$. Thus, $G_1' = wt$ where $w$ is adjacent to $x$ and $t$ is adjacent to $z$, and $G_2' = ab$ where $a$ is adjacent to $z$ and $b$ is adjacent to $y$. However, $G$ contains a maximal matching of size $3$ and yet $\{wt, wx, ab, by\}$ is a minimal edge dominating set, which is a contradiction. 

Thus, we may assume that $G_1'$ is the only component in $G'$, and from previous arguments, $G_1' \in \{K_2, C_4\}$. However, $G_1' \ne K_2$ for this would imply that $\deg_G(y) = 2$. Thus, $G_1' = C_4$ and  one can easily verify that $G = Cr$. 

\end{proof}

Using all of the previous results in this section, we are now ready to prove Theorem~\ref{thm:onetriangleWED}. 
\vskip5mm

\noindent \textbf{Theorem~\ref{thm:onetriangleWED}.} \emph{
$G$ is a connected, well-edge-dominated graph with exactly one triangle if and only if $G \in \mathcal{T}\cup \mathcal{F} \cup \{K_3, Cr, \mathcal{H}, \mathcal{DH}\}$.
}
\vskip2mm
\begin{proof}
We shall assume that $G \ne K_3$ is well-edge-dominated. Let $xyz$ be the lone triangle in $G$. If exactly one vertex on the triangle has degree $2$ in $G$, then by Theorem~\ref{thm:degree2} $G \in \{\mathcal{H}, \mathcal{DH}\}$. If exactly two vertices on the triangle have degree $2$ in $G$, then without loss of generality, we may assume that $\deg_G(x) = \deg_G(y) = 2$. Consider $G' = G - \{x, y, z\}$. By Lemma 4, $G'$ contains nontrivial components. Thus, by Theorem~\ref{thm:houseonfire}, $G \in \mathcal{T}$. Therefore, we may assume that each of $x, y$, and $z$ have degree at least $3$ in $G$. If $G - N_e[f]$ contains at least two nontrivial components for any edge $f$ on $xyz$, then $G \in \mathcal{F}$ by Theorem~\ref{thm:houseonfire2}. So we may assume that there is some edge $f$ on $xyz$ such that $G-N_e[f]$ contains only one nontrivial component. Relabeling if necessary, we shall assume that $G-N_e[xy]$ contains one nontrivial component. By Theorem~\ref{thm:crystal}, $G = Cr$.
\end{proof}

\section{Outerplanar Graphs}\label{sec:OP}
 To show that there are well-edge-dominated graphs containing multiple triangles, we now focus on outerplanar graphs. Recall that a graph $G$ is \emph{outerplanar} if $G$ has a planar drawing with every vertex on the outer face. A \emph{fan} on $n$ vertices, denoted $F_n$, is obtained by taking the join of $K_1$ and $P_{n-1}$. 

\begin{theorem}\label{thm:equiouter}
If $G$ is outerplanar with $|V(G)| = 2k$ for $k\geq3$ or $|V(G)| = 2m+1$ for $m\geq4$, then $G$ is not equimatchable. 
\end{theorem}

 \begin{proof}
Label the vertices of $G$ as $v_1v_2\dots v_{n(G)}$ where each $v_iv_{i+1}$ (taken modulo $n(G)$) is an edge on the outer face of $G$ and $\deg_G(v_1) = 2$. When $|V(G)| = 2k$, $\{v_iv_{i+1} : i \in [2k], \text{ $i$ odd}\}$ is a perfect matching in $G$. Alternatively, consider the matching $M_2 = \{v_2v_3\} \cup \{v_iv_{i+1}: 5 \le i < 2k, \text{ $i$ odd}\}$. Note that $M_2$ is maximal since the only vertices not saturated by $M_2$ are $v_1$ and $v_4$ yet there is no edge edge $v_1$ and $v_4$ since $N_G(v_1) = \{v_2, v_{2k}\}$. Thus, $G$ is not equimatchable.

Hence, we consider when  $|V(G)| = 2m+1$ for $m\geq4$.  Suppose first that $\{v_1, v_4, v_7\}$ is independent in $G$. Then $M_1 = \{v_iv_{i+1}: 2 \le i \le 2m, \text{ $i$ even}\}$ saturates every vertex except $v_1$. Therefore, $M_1$ is a maximal matching. Alternatively, consider $M_2 = (M_1 - \{v_4v_5, v_6v_7\}) \cup \{v_5v_6\}$. $M_2$ dominates every edge except possibly edges between vertices in $\{v_1, v_4,v_7\}$. Since this is an independent set, $M_2$ is in fact a maximal matching of cardinality $|M_1| - 1$. Thus, $G$ is not equimatchable. 

So we shall assume that  $\{v_1, v_4, v_7\}$ is not an independent set. It follows that  $v_4v_7 \in E(G)$ as we have assumed $N_G(v_1) = \{v_2, v_{2m+1}\}$. Suppose first that $|V(G)| \ge 11$. If $\{v_2, v_5, v_8\}$ is an independent set, both $M_1 = \{v_iv_{i+1} : 3 \le i \le 2m+1, \text{ $i$ odd}\}$ and $M_2 = (M_1 - \{v_5v_6, v_7v_8\}) \cup \{v_6v_7\}$ are maximal matchings implying $G$ is not equimatchable. Thus, the only case left to consider is when $\{v_2, v_5, v_8\}$ is not an independent set. In this case, we know that $v_5$ is not adjacent to $v_2$ or $v_8$ as we began with a plane drawing where $v_1 \dots v_{2m+1}$ is the outer boundary. Therefore, $v_2v_8 \in E(G)$ implying that $\{v_3, v_6, v_9\}$ is an independent set. In this case, $N_G(v_3) \subseteq \{v_2, v_4, v_7, v_8\}$ and $M_1 = \{v_iv_{i+1}: 4 \le i \le 2m, i \text{ even}\}\cup \{v_1v_2\}$ is a maximal matching in $G$. Alternatively, $M_2 = (M_1 - \{v_6v_7, v_8v_9\}) \cup \{v_7v_8\}$ is a maximal matching as the only vertices which are not saturated are $\{v_3, v_6, v_9\}$. It follows that $G$ is not equimatchable if $|V(G)| \ge 11$. 

The only case left to consider is when $|V(G)| = 9$ and $\deg_G(v_1) = 2$ and $v_4v_7 \in E(G)$. As above, if $\{v_2, v_5, v_8\}$ is an independent set in $G$, then $\{v_3v_4, v_5v_6, v_7v_8, v_1v_9\}$ and $\{v_3v_4, v_6v_7, v_1v_9\}$ are both maximal matchings. So we shall assume that $v_2v_8 \in E(G)$. It follows that $\{v_3, v_6, v_9\}$ is an independent set and both $\{v_1v_9, v_2v_8, v_4v_7, v_5v_6\}$ and $\{v_1v_2, v_4v_5, v_7v_8\}$ are maximal matchings in $G$.

 \end{proof}

\begin{corollary} Let $G$ be outerplanar.  $G$ is well-edge dominated if and only if $G \in \{C_3, C_4, C_5, \mathcal{H}, F_5, C_7, C_7^*, \mathcal{DH}\}$. 
\end{corollary}

\begin{proof}
 By Theorem~\ref{thm:equiouter}, if $G$ is well-edge-dominated, then $|V(G)| \in \{3, 4, 5, 7\}$. Computer search yields all well-edge-dominated outerplanar graphs. 
\end{proof}

\section{Acknowledgements}
We would like to thank Yale's SUMRY program for supporting this project. C. Kaneshiro was supported by Princeton University’s Office of Undergraduate Research OURSIP Program through the Hewlett Foundation.

\end{document}